\documentclass[leqno, 12pt]{article}

\usepackage{amssymb, amsmath, amsthm} 

 \addtocounter{section}{-1}

\newcommand{\on}{\operatorname}
\newcommand{\mc}{\mathcal}
\newcommand{\mf}{\mathfrak}

\newcommand{\be}{\begin{equation}}
\newcommand{\ee}{\end{equation}}

\def\lra{{{\,\,\longrightarrow\,\,}}} 
 
\newcommand{\ul}{\underline}

\newcommand{\CC}{{\mathbb C}} 
\newcommand{\DD}{{\mathbb D}}

\newcommand{\GG}{{\mathbb G}}

\newcommand{\PP}{{\mathbb P}} 
 
\newcommand{\RR}{{\mathbb R}} 

\newcommand{\Ac}{{\mc A}}

\newcommand{\Ec}{{\mc E}}

\newcommand{\Mc}{{\mc M}} 
 
\newcommand{\Oc}{{\mc O}} 
\newcommand{\Pc}{{\mc P}} 
\newcommand{\Qc}{{\mc Q}} 
 
\newcommand{\Sc}{{\mc S}} 
 
\newcommand{\Uc}{{\mc U}}

\newcommand{\Uen}{{\mf U}}

\newcommand{\cO}{{\Oc}}
\newcommand{\cP}{{\Pc}}

\newcommand{\cA}{{\Ac}}

\newcommand{\cE}{{\Ec}}

\newcommand{\cQ}{{\Qc}}

\newcommand{\Aut}{\operatorname{Aut}}
\newcommand{\End}{\on{End}}

\newcommand{\MHS}{\operatorname{MHS}}
\newcommand{\gr}{\operatorname{gr}}
\newcommand{\barF}{\overline{F}\, }
\newcommand{\Bun}{\operatorname{Bun}}
\newcommand{\ds}{/ \hskip -3pt /}
\newcommand{\Vect}{\operatorname{Vect}}
\newcommand{\Ree}{\operatorname{Re}}
\newcommand{\Hod}{\operatorname{Hod}}
\newcommand{\Hom}{\operatorname{Hom}}

\newcommand{\FLie}{\operatorname{FLie}}
\newcommand{\Lie}{\operatorname{Lie}}
\newcommand{\cl}{\operatorname{cl}}
\newcommand{\pol}{\operatorname{pol}}
\newcommand{\ad}{\operatorname{ad}}
\newcommand{\Rep}{\operatorname{Rep}}
\newcommand{\Spec}{\operatorname{Spec}}
\newcommand{\Proj}{\operatorname{Proj}}
\newcommand{\Rs}{\operatorname{Rs}}
\newcommand{\PR}{{\Bbb P}\operatorname{Rs}}

 \author{Mikhail Kapranov\thanks{mikhail.kapranov@yale.edu} }

\title{Real mixed Hodge structures.}

\begin{document}

\maketitle

\newtheorem{thm}[equation]{Theorem}
\newtheorem{cor}[equation]{Corollary}
\newtheorem{lem}[equation]{Lemma}
\newtheorem{prop}[equation]{Proposition}
\newtheorem{refo}[equation]{Reformulation}

\newtheoremstyle{example}{\topsep}{\topsep}%
     {}
     {}
     {\bfseries}
     {.}
     {2pt}
     {\thmname{#1}\thmnumber{ #2}\thmnote{ #3}}

   \theoremstyle{example}
   \newtheorem{nota}[equation]{Notation}
   \newtheorem{Defi}[equation]{Definition}
   \newtheorem{rem}[equation]{Remark}
   \newtheorem{rems}[equation]{Remarks}
   \newtheorem{prob}[equation]{Problem}
   \newtheorem{exas}[equation]{Examples}
   \newtheorem{exa}[equation]{Example}
  \numberwithin{equation}{section}

\section { Introduction.}

\vskip 1cm

Let $\MHS_{\mathbb R}$ be the category of real mixed Hodge structures 
\cite{D1}\cite{D2}. Thus, an object of $\MHS_{\Bbb R}$ consists of a finite-dimensional
$\Bbb R$-vector space $V$, equipped with an increasing filtration
$W_\bullet$ and a decreasing filtration $F^\bullet$
on $V_{\Bbb C} = V\otimes_{\Bbb R} \Bbb C$. These filtrations are required
to satisfy the following condition:
\be\label{mixed-hodge-equation}
\gr^p_F \gr^q_{\overline{F}} \gr_n^W \, (V_{\Bbb C}) = 0
\quad {\operatorname{ for}} \quad n\neq p+q.
\ee
Here $\barF^\bullet$ is the filtration obtained from $F^\bullet$ by 
complex conjugation, and we denote the filtration induced by $W_\bullet$
on $V_{\Bbb C}$ by the same symbol $W_\bullet$. 

\vskip .3cm

It is known that $\MHS_{\Bbb R}$ is an abelian tensor category.
The goal of this note is to give a ``gauge-theoretic''
 description of $\MHS_{\Bbb R}$. 

Consider $\Bbb C$, the set of complex numbers, as a 2-dimensional
real algebraic variety (i.e., an algebraic variety over $\RR$),
 and the multiplicative group $\Bbb C^*$ as
a 2-dimensional real algebraic group acting on this variety. In more
formal algebro-geometric terms, we are applying the Weil restriction
functor $R_{\Bbb C /\Bbb R}$ from $\Bbb C$ to $\Bbb R$:
\be 
{\Bbb C} = R_{\Bbb C /\Bbb R} ({\Bbb A}^1) = {\Bbb A}^2_{\Bbb R}, \quad 
{\Bbb C}^* = 
S := R_{\Bbb C /\Bbb R}({\Bbb G}_m), 
\ee
see \cite{D1}  (2.1.2). Let $\Bun_\nabla({\Bbb C}; {\Bbb C}^*)$ be the category of
(real) algebraic vector bundles on $\Bbb C$, equivariant with
respect to the ${\Bbb C}^*$-action and equipped with an invariant
(real algebraic) connection. Note that the connections are not assumed flat.


\begin{thm}\label{main-theorem}  There exists an equivalence of abelian tensor
categories
$$h: \Bun_\nabla({\Bbb C}; {\Bbb C}^*)\lra \MHS_{\Bbb R}$$
 with the following properties:

(a)  If $(E,\nabla)$ is an object of $\Bun_\nabla({\Bbb C}; {\Bbb C}^*)$
and $(V, W_\bullet, F^\bullet,
\barF^\bullet) = h(E,\nabla)$,  then the space $\gr_\bullet^W (V)$ is canonically identified with $E_0$, the
fiber of $E$ at $0$.

(b) In the situation of (a), the  vector space $V$ is identified with
the space of covariantly constant sections
$$V\quad  =\quad  H^0_\nabla\bigl( \bigl\{ \Ree (t) = -1/2\bigr\}, E\bigr).$$

(c)  $\nabla$ is flat if and only if $V$ is split, i.e., isomorphic to
a direct sum of pure Hodge structures. 

(d)  The functor of ``absolute cohomology''
$$R\Gamma_{\Hod}: \MHS \to \Vect_{\Bbb R}, \quad V\mapsto  R\Hom_{\MHS}
({\Bbb R}(0), V)$$
of Beilinson \cite{B} is identified  with the invariant part of the  ``de Rham sequence'' functor
$$(E,\nabla) \quad  \mapsto\quad \Gamma\bigl( {\Bbb C}, \,  \bigl\{ E\buildrel\nabla\over\longrightarrow
 E\otimes\Omega^1_{\Bbb C}\bigr\}\bigr)^{{\Bbb C}^*}.$$
Here $\Omega^1_{\Bbb C}$ is the sheaf of real 1-forms on $\Bbb C$ as a real 2-manifold.
\end{thm}

\begin{rems} (a) This gives a geometric
interpretation of the pro-unipotent group  $\frak U$ which, according to Deligne
 \cite{D2}, 
governs mixed Hodge structures. Indeed, $\frak U$
is realized as (the pro-algebraic completion of) the group
of piecewise smooth loops in $\Bbb C$  considered up to reparametrization and cancellation
(with the group operation being the composition of loops). This group
of loops acts in any bundle with connection via the holonomy. 
See Subsection \ref{fundamental-algebra-subsection}
 below for a discussion at the Lie algebra level. 

\vskip .1cm

(b)  Theorem \ref{main-theorem} bears interesting  similarities with the results of
Connes and Marcolli \cite{Connes-Marcolli}. They considered a certain
tensor category of ``equisingular" families of connections and identified 
it with the category
of representations of a
pro-unipotent group $\Uc^*$. This group contains $\GG_m$ and   the commutant  
$[\Uc^*, \Uc^*]$ is free
on generators of $\GG_m$-weights 1,2,3,... On the other hand, our results
show that the category $\Bun_\nabla(\CC, \CC^*)$ is equivalent to the
category of representations of the group $\frak U$ of Deligne. This group contains
$\CC^*=S$, and $[\Uen,\Uen]$ is free on generators $z_{p,q}$, $p,q\geq 1$, 
 of bi-weight $p,q$. So the relation between $\Uen$ and $\Uc^*$ is
 qualitatively the same as the relation between all mixed Hodge structures
 and iterated extensions of Tate structures. However, the two
 categories of connections leading to these groups are quite different
 and a direct relation of \cite{Connes-Marcolli}
 with Hodge structures is not known.

\vskip .1cm

(c) The group $\CC^*$ does not act freely on $\CC$ because of the fixed point
at 0, so from the topological point of view the quotient $\CC/\CC^*$ is ill behaved
and should be replaced by the {\em quotient stack} $\CC\ds\CC^*$. 
This stack can be seen 
 as
a local archimedeal analog of the ``adele class space''
$A\ds k^*$ of Connes \cite{C}.
 Here $k$ is a global field, and $A$ is
its ring of adeles. Quotients similar to ${\Bbb C}\ds{\Bbb C}^*$ were
also considered by Laumon \cite{L} from the point of view
of D-modules and constructible sheaves. Categories of
equivariant connections on other pre-homogeneous vector spaces
can provide interesting analogs of the category of mixed Hodge structures.

\end{rems}



\vskip .3cm

The proof  of Theorem \ref{main-theorem} will be given in \S \ref{construction-section}.
 It is based on a version of the 
Radon-Penrose transform
which takes equivariant bundles with connections on (the complexification of) $\Bbb C$
into equivariant  algebraic vector  bundles on a punctured projective plane. The latter
bundles can be realized as Rees bundles of mixed Hodge structures, following
the work of Penacchio \cite{P1}\cite{P2}. 

 In \S \ref{Deligne-operator-section} we will relate our description
with the description of Deligne \cite{D2} using the ``Hodge monodromy operator''.
We also give an interpretation of Deligne's Hodge group $\frak G$ and its
Lie algebra. 

\vskip .3cm

I am deeply grateful to A. B. Goncharov who explained to me some
basics of mixed Hodge structures. He also informed me about
his work  \cite{G} (then in progress). Theorem \ref{main-theorem} turned out to be
very closely related to the approach of \cite{G}. In particular,
the ``twistor line'' of \cite{G} is naturally identified with the real line $\Ree(t)=-1/2$
from Theorem \ref{main-theorem}(c). See Remark \ref{twistor-line-remark}
below for more details. 
I would also like to thank F. Loeser for pointing out the
work \cite{L}. This work was partially supported by an NSF grant.

  \numberwithin{equation}{subsection}

\section {Construction of the equivalence.}\label{construction-section}

\subsection {Complex version.}\label{complex-version-subsection}
We will deduce Theorem \ref{main-theorem} from its complex version.
Recall that a {\em complex mixed Hodge structure} is 
a datum consisting of a finite-dimensional complex vector space $V$ with
an increasing filtration $W$ and two decreasing filtrations $F', F''$
which satisfy the condition
\be \gr^p_{F'} \gr^q_{{F''}} \gr_n^W \, (V) = 0
\quad {\operatorname{ for}} \quad n\neq p+q.
\ee
similar to \eqref{mixed-hodge-equation}.
  Another name for such an object
is ``a triple of opposite filtrations'' \cite{D1}. We denote
by $\MHS_\CC$ the category of complex mixed Hodge structures.
Like $\MHS_\RR$, it is an abelian tensor category. 

\vskip .2cm

For a complex vector space $V$ we denote $\sigma(V)$ the complex conjugate space.
It consists of symbols $\sigma(v), v\in V$ with the operations
\be 
\sigma(v)+\sigma(v')=\sigma(v+v'), \quad \lambda\sigma(v) = \sigma(\overline{\lambda}v), \, 
\lambda\in {\Bbb C}. 
\ee
Thus we have an antilinear isomorphism
\be\label{sigma-map-equation}
\sigma: V\to\sigma(V). 
\ee
Given a  complex mixed Hodge structure $(V, W_\bullet, F'{}^\bullet, F''{}^\bullet)$, its
{\em complex conjugate structure}  is defined by 
\be\label{conjugate-structure-equation}
\sigma(V, W_\bullet, F'{}^\bullet, F''{}^\bullet) = (\sigma(V), \sigma(W_\bullet),
 \sigma(F''{}^\bullet),
\sigma( F'{}^\bullet)).
\ee
This defines an action of the  Galois group
\be 
\Gamma = \operatorname{Gal} ({\Bbb C}/{\Bbb R}) = \{ 1,\sigma\}
\ee
 on the category $\MHS_{\Bbb C}$, 
 and $\MHS_{\Bbb R}$ consists of $\Gamma$-equivariant
objects in $\MHS_{\Bbb C}$. So we describe $\MHS_{\Bbb C}$ in a way
compatible with the $\Gamma$-action. 

\vskip .2cm

Consider the affine plane ${\Bbb A}^2 = \operatorname{Spec} 
\,{\Bbb C}[t_1, t_2]$ over $\Bbb C$ 
with the standard action of the torus ${\Bbb G}_m^2$. 
Let $\Bun_\nabla({\Bbb A}^2; {\Bbb G}_m^2)$ be the category of
complex algebraic vector bundles  on ${\Bbb A}^2$ which
are equivariant with respect to the ${\Bbb G}_m^2$-action
and equipped with an equivariant (complex algebraic) connection.
There is a $\Gamma$-action on $\Bun_\nabla({\Bbb A}^2; {\Bbb G}_m^2)$
induced by the action on ${\Bbb A}^2$ given by
\be\label{action-A2-equation}
\sigma(t_1, t_2) = (\overline{t}_2, \overline{t}_1), 
\ee
and by a similar action on ${\Bbb G}_m^2$. The following then implies Theorem
 \ref{main-theorem}(a). 

\begin{thm}\label{main-complex-theorem}
 $\MHS_{\Bbb C}$ is equivalent to 
$\Bun_\nabla({\Bbb A}^2; {\Bbb G}_m^2)$ as a tensor category with
$\Gamma$-action. 
\end{thm}

The proof will be given in Subsection \ref{main-theorems-proof-subsection}
below.

\subsection{ The Radon-Penrose transform of equivariant
 connections.}\label{radon-penrose-subsection}

Let us compactify
 ${\Bbb A}^2 = \Spec\, {\Bbb C}[t_1, t_2]$
  to 
  ${\Bbb P}^2 = \Proj\, {\Bbb C}[u_0, u_1, u_2]$,
so ${\Bbb A}^2$ is given by $u_0\neq 0$, and
\be 
{\Bbb P}^2 \,\,=\,\, {\Bbb A}^2 \cup {\Bbb P}^1_\infty, \quad {\Bbb P}^1_\infty\,\, =\,\, \Proj\, 
{\Bbb C}[t_1, t_2],\quad t_i = u_i/u_0.
\ee
Let $\check{\Bbb P}^2$ be the dual projective plane of lines in ${\Bbb P}^2$,
so $\check{\Bbb P}^2 = \Proj\, {\Bbb C}[v_0, v_1, v_2]$, where $(v_i)$ are the dual coordinates.
As usual, points in ${\Bbb P}^2$ give lines in $\check{\Bbb P}^2$: for a point $x\in {\Bbb P}^2$
we denote $\lambda_x$ the set of lines in ${\Bbb P}^2$ through $x$. 
Now, lines in ${\Bbb A}^2$ form
$$\check{\Bbb P}^2_0\,\,\, =\,\,\, \check {\Bbb P}^2 - \{[1:0:0]\},$$
so we have the incidence diagrams
\be\label{incidence-diagram-equation}
\begin{matrix} {\Bbb A}^2&\buildrel q_0\over\longleftarrow&Q_0&\buildrel p_0\over\longrightarrow&
\check{\Bbb P}^2_0\\
\big\downarrow k&&\big\downarrow&&\big\downarrow j\\
{\Bbb P}^2&\buildrel q\over\longleftarrow &Q&\buildrel p\over\longrightarrow &\check{\Bbb P}^2
\end{matrix}
\ee

\begin{prop}\label{algebraic-solutions-proposition}
 Let  $(E,\nabla)\in\Bun_\nabla({\Bbb A}^2;  {\Bbb G}_m^2)$ and $r=\on{rk}(E)$.
 Consider
the  following sheaf of $\cO$-modules on $\check{\Bbb P}^2_0$:
$$\cE^0 = p_{0*}^\nabla(q_0^*E),$$
where $p_{0*}^\nabla$ is the subsheaf in the direct image $p_{0*}$ consisting of sections
covariantly constant along the fibers of $p_0$. Then $\cE^0$ is locally free
of rank $r$, so it is an algebraic vector bundle on $\check{\Bbb P}^2_0$. 
\end{prop}

 We call $\cE^0$ 
the {\it Radon-Penrose transform} of $(E,\nabla)$, following \cite{M}, Ch. 2, \S 2. 

\vskip .3cm

\noindent {\sl Proof of the proposition:} The variety $\check{\Bbb P}^2_0$ is covered
by two open charts isomorphic to $\mathbb A^2$. One of them, $U=\Spec\, \CC[a,b]$,
parametrizes lines of the form 
\be\label{at+b-equation}
t_2=at_1+b,
\ee
 and the other one is defined similarly,
with the roles of $t_1, t_2$ exchanged. So we will prove that $\Ec^0|_U$ is
a free $\Oc_U$-module of rank $r$. 

\vskip .2cm

Let $V$ be the space of $\GG_m^2$-invariant sections of $E$ over $\GG_m^2$,
so $\dim(V)=r$. Over $\GG_m^2$, we have a trivialization $E\simeq \Oc_{\GG_m^2}\otimes V$.  
With respect to this trivialization, the equivariant connection $\nabla$ has the form
\be\label{nabla-equivariant-equation}
\nabla\,\,=\,\,d + B_1 \,d\log t_1 + B_2 \,d\log t_2,\quad B_1, B_2\in\End(V).
\ee
Now, Eq.  \eqref{at+b-equation} identifies $p_0^{-1}(U)$ with $U\times\mathbb A^1$,
where
$\mathbb A^1=\Spec\, \CC[t_1]$. Denote by $\widetilde E$ the algebraic vector bundle
on $U\times\mathbb A^1$ corresponding to the bundle $q_0^{*}(E)|_{p_0^{-1}(U)}$ under this
identification. Then 
 $\nabla$ induces a  relative connection in $\widetilde E$
 along the fibers of the projection $U\times\mathbb A^1\to U$. 
 Denote this relative connection by $D$.
 Substituting \eqref{at+b-equation}
into \eqref{nabla-equivariant-equation} and using the above trivialization of
$E$ on $\GG_m^2$, we find a trivialization of $\widetilde E$
near $U\times\{\infty\}\subset U\times\PP^1$ such that the connection matrix of $D$
in this trivialization
has at most first order pole near $t_1=\infty$.

This means that for each choice of numerical values $a,b\in\CC$ we have then a nonsingular
connection in an algebraic vector bundle  
$\widetilde E_{a,b}=E|_{\{t_2=at_1+b\}}$
on $\mathbb A^1$, having a regular singularity at infinity, see
\cite{Deligne-sing}, Th. 1.1.2(i). As such, it has 
 a fundamental solution whose matrix elements
are regular functions (polynomials) in $t_1$, see \cite{Deligne-sing},
Th. 1.1.9. As $a,b$ vary, the coefficients
 of these polynomials are
regular functions in $a,b$ as they are found by the standard recursive formulas. 
This implies that 
the relative connection $D$ in $\widetilde E$ is algebraically trivial
(isomorphic to a pullback of a bundle on $U$, with trivial
relative connection), and
so the sheaf of covariantly constant sections is a free $\Oc_U$-module
of rank $r$. \qed

\vskip .3cm

The action of the torus ${\Bbb G}_m^2$ on ${\Bbb A}^2$ by dilations extends to ${\Bbb P}^2$
and, by duality, to $\check{\Bbb P}^2$. Both ${\Bbb P}^2$ and $\check{\Bbb P}^2$ are toric
varieties with respect to this action, and $\check{\Bbb P}^2_0$ is a torus invariant open set. 
By construction, the vector bundle $\cE^0$ on $\check{\Bbb P}^2_0$ is ${\Bbb G}_m^2$-equivariant. 

\vskip .2cm

Note further that the embedding $j: \check{\Bbb P}^2_0\to \check{\Bbb P}^2$ misses just one
point, so the direct image
\be\label{extension-E-equation}
 \cE = j_* \cE^0
\ee
is a reflexive coherent sheaf on the  surface $\check{\Bbb P}^2$ and thus a vector bundle.
This bundle is still equivariant with respect to ${\Bbb G}_m^2$. 

We denote by $\cQ$ the category 
of vector bundles $\cE$ on $\check{\Bbb P}^2$
equivariant under ${\Bbb G}_m^2$ and trivial on each line $\lambda\subset {\Bbb P}^2$
except perhaps the lines corresponding to points of ${\Bbb P}^1_\infty\subset {\Bbb P}^2$.
We have a natural structure of a tensor category on $\cQ$ given by
tensor product of vector bundles.

\begin{thm}\label{radon-penrose-theorem}
The correspondence $(E,\nabla)\mapsto \cE$
establishes an equivalence of tensor categories
 $$\Bun_\nabla({\Bbb A}^2; {\Bbb G}_m^2)\lra \cQ.$$
\end{thm}

\noindent {\sl Proof:} We already showed how to construct $\cE$ from $(E,\nabla)$ and
that $\cE$ is equivariant. Note that $\cQ$ is equivalent, via 
\eqref{extension-E-equation}, to the category $\cQ_0$
of equivariant bundles on $\check{\Bbb P}^2_0$ trivial on {\em all the lines contained
in} $\check {\Bbb P}^2_0$. Indeed, these are precisely the lines not 
corresponding to the points of ${\Bbb P}^1_\infty$. 
 Now, the equivalence of $\Bun_\nabla({\Bbb A}^2; {\Bbb G}_m^2)$
with $\cP_0$ is an equivariant version of a particular case of the general fact about 
Radon-Penrose transforms (\cite{M}, Ch. 2, \S 2, Theorem 2.3). Indeed, the first line in
\eqref{incidence-diagram-equation} is a particular case of a double fibration considered
 in \cite{M}, Ch. 2, Sect. 2.1.
The cited theorem establishes an equivalence between holomorphic bundles on 
$\check {\Bbb P}^2_0$
trivial on
all the lines and holomorphic bundles with connections on ${\Bbb A}^2$ but in
 our equivariant situation
we can restrict to algebraic bundles on both sides, in virtue of Proposition
\ref{algebraic-solutions-proposition}. 
 Let us just explain why
$\cE$ is trivial on each line contained in $\check{\Bbb P}^2_0$, i.e., any line $\lambda$ of
the form $\lambda_x$, $x\in {\Bbb A}^2$. Indeed, if $l\subset {\Bbb A}^2$
is a line through $x$, then the restriction  of covariantly constant sections to $x$
gives an isomorphism
$$\cE_l = H^0_\nabla(l, E)\buildrel \simeq\over\longrightarrow E_x,$$
so
\be\label{E-lambda-equation}
\cE|_\lambda = \cO_\lambda\otimes E_x
\ee
is trivial. \qed

\subsection  {Rees bundles and the work of Penacchio.} 
We now relate Theorem \ref{radon-penrose-theorem}
 with the description of $\MHS_{\Bbb C}$ given by
Penacchio \cite{P1}\cite{P2}. 
Denote by
\be\label{rees-bundle-equation}
 \begin{aligned}
\Rs(V, W_\bullet, F'{}^\bullet, F''{}^\bullet)\quad  =\quad  \bigoplus_{i,j,k\in {\Bbb Z}}
(W_i\cap F'{}^{-j}\cap F''{}^{-k})(V) v_0^i v_1^j v_2^k \\
\quad \subset\quad  
{\Bbb C}[v_0^{\pm 1}, v_1^{\pm 1}, v_2^{\pm 1}]
\end{aligned}
\ee
the Rees module over ${\Bbb C}[v_0, v_1, v_2]$ corresponding to the 3-graded
vector space $(V, W_\bullet, F'{}^\bullet, F''{}^\bullet)$. Here the minus signs
before $j$ and $k$ correpond to converting the decreasing fultrations
$F', F''$ into increasing ones. Both ${\Bbb C}[v_0, v_1, v_2]$ and
$\Rs(V, W_\bullet, F'{}^\bullet, F''{}^\bullet)$ have compatible ${\Bbb Z}^3$-gradings
which translate to a ${\Bbb G}_m^3$-action on $\check{\Bbb A}^3 = \Spec \, {\Bbb C}[v_0, v_1, v_2]$
and into an equivariant structure of the coherent sheaf on $\check{\Bbb A}^3$
corresponding to $\Rs(V, W_\bullet, F'{}^\bullet, F''{}^\bullet)$. If we consider
the gradings by total degree, we can form the projective plane $\check{\Bbb P}^2 = 
\Proj\, {\Bbb C}[v_0, v_1, v_2]$ which we identify with the $\check{\Bbb P}^2$
from Subsection \ref{radon-penrose-subsection},
 and a coherent sheaf on $\check{\Bbb P}^2$ corresponding to the
graded module $\Rs(V, W_\bullet, F'{}^\bullet, F''{}^\bullet)$. We denote this
sheaf $\PR(V, W_\bullet, F'{}^\bullet, F''{}^\bullet)$ or simply by $\cE$. The quotient torus
\be 
{\Bbb G}_m^2 = {\Bbb G}_m^3/{\Bbb G}_m
 \ee
(quotient by the diagonal embedding) is identified with the ${\Bbb G}_m^2$ acting
on $\check{\Bbb P}^2$, and $\cE$ is equivariant. Further, it is known that  Rees
modules are reflexive, so $\cE$ is a ${\Bbb G}_m^2$-equivariant vector bundle on $\check {\Bbb P}^2$. 
Since the variables $v_0, v_1, v_2$ are associated to the filtrations $W, F', F''$ in 
\eqref{rees-bundle-equation},
we denote the coordinate lines in $\check{\Bbb P}^2$ by
\be 
\check{\Bbb P}^1_W = \{ v_0=0\}, \quad \check{\Bbb P}^1_{F'} = \{ v_1=0\}, \quad
\check{\Bbb P}^1_{F''}=\{ v_2=0\},
\ee
and the torus fixed points by
\be 
\check{\Bbb P}^0_{WF'}= \{ v_0=v_1=0\}, \quad \check{\Bbb P}^0_{W F''} = \{ v_0=v_2=0\},
\quad
\check {\Bbb P}^1_{F'F''}=\{v_1=v_2=0\}.
 \ee
Then, the restrictions of the Rees bundles to the coordinate lines are found as follows
(see \cite{P2}, (2.6.2)):
\be\label{rees-restriction-equation} 
\cE|_{\check{\Bbb P}^1_W} = \bigoplus_n \PR (\gr_n^WV, F'{}^\bullet, F''{}^\bullet)\otimes
\cO(-n),
\ee
where on the RHS we have the Rees bundles on ${\Bbb P}^1$ corresponding to the bifiltered spaces
$(\gr_n^W V, F'{}^\bullet, F''{}^\bullet)$ (induced filtrations). Similarly for other
coordinate lines. The restrictions to (i.e., fibers over) the fixed points are given by the bigraded spaces associated
to the corresponding  pairs of filtrations:
\be 
\cE|_{\check{\Bbb P}^0_{WF'}} = \gr_\bullet^W \gr^\bullet_{F'}(V), \quad {\operatorname {etc.}}
 \ee
Now, recall an observation of Deligne, used by Simpson \cite{S}.

 \begin{prop}
Let $F'{}^\bullet, F''{}^\bullet$ be two filtrations on a finite-dimensional
vector space $V$. Then the following are equivalent:
\hfill\break
(i) $F'{}^\bullet$ and $F''{}^\bullet$ are $n$-opposite (induce a pure complex Hodge structure of weight $n$).
\hfill\break
(ii) The Rees bundle $\PR(V, F'{}^\bullet, F''{}^\bullet)$ on $P^1$ is a direct
sum of several copies of $\cO(n)$.
\end{prop}

As observed in \cite{P2}, the isomorphism 
\eqref{rees-restriction-equation} implies that for a complex mixed
Hodge structure $(V, W_\bullet, F'{}^\bullet, F''{}^\bullet)$ the restriction
$\cE|_{\check {\Bbb P}^1_W}$ is trivial:
\be\label{E-gr-equation} 
\cE|_{\check {\Bbb P}^1_W} \simeq \gr_\bullet^W(V)\otimes \cO_{\check {\Bbb P}^1_W}.
 \ee

Denote by $\cP$ the following category. 
 Objects of $\cP$ are ${\Bbb G}_m^2$-equivariant
vector bundles $\cE$ on ${\Bbb P}^2$ which are trivial on $\check{\Bbb P}^1_W$.
Morphisms in $\cP$ are equivariant morphisms of bundles which have constant rank
everywhere except possibly the point $\check{\Bbb P}^0_{F' F''}$. 
Clearly, $\cP$ is a tensor category with respect to the tensor product of
vector bundles. 
The main result of \cite{P1}\cite{P2} is (see \cite{P1},  Th. 3.1):

 \begin{thm}\label{rees-bundle-theorem}
  The Rees bundle construction defines an equivalence of
tensor categories
$$\MHS_{\Bbb C}\to \cP.$$
 \end{thm}

Notice now the following:

 \begin{lem}\label{P-Q-lemma}
  The categories $\cP$ and $\cQ$ are equivalent. 
\end{lem}

\noindent {\sl Proof:} We first identify the objects. On the surface of it, 
$\cP$ seems to have more objects, as we require triviality on one line only,
rather than on all lines not meeting $\check{\Bbb P}^0_{F'F''}$.
However, triviality is an open condition for vector bundles  on ${\Bbb P}^1$.
So if $\cE\in \cP $ is trivial on $\check{\Bbb P}^1_W$, it is trivial on
an open subset $U$ of lines in $\check{\Bbb P}^2$ containing $\check{\Bbb P}^1_W$.
By equivariance, $U$ can be assumed to  be preserved under the torus action. This implies
that $U$ contains all lines not meeting $\check{\Bbb P}^0_{F'F''}$.
So $\cE$ is an object of $\cQ$ as well.

\vskip .1cm

Next, we identify the  morphisms. On the surface of it, $\cQ$ seems to have more morphisms, as
we do not require the constant rank condition. So let $f: \cE\to \cE'$ be a morphism
in $\cQ$, i.e., just an invariant morphism of equivariant vector bundles, both being
objects of $\cQ$. Let $\lambda\subset
\check{{\Bbb P}^2}$ be a projective line 
such that both $\cE$ and $\cE'$ are trivial on $\lambda$.
Then clearly $f$ has constant rank along $\lambda$. On the other hand, 
 any two points of $\check{\Bbb P}^2 -
\check{\Bbb P}^0_{F' F''}$ are connected by a chain of lines $\lambda$ as above.
This implies that the rank of $f$ is constant on $\check{\Bbb P}^2 -
\check{\Bbb P}^0_{F'F''}$, so $f$ is a morphism of $\cP$. \qed

 \subsection{Proof of Theorems \ref{main-complex-theorem}  and 
 \ref{main-theorem}(a)-(d).}\label{main-theorems-proof-subsection}
Combining now Theorems \ref{radon-penrose-theorem}, \ref{rees-bundle-theorem}
 and Lemma \ref{P-Q-lemma}, we get an equivalence
of tensor categories
\be \MHS_{\Bbb C}\buildrel \PR\over\longrightarrow \cQ = \cP \buildrel \Psi\over\longrightarrow 
\Bun_\nabla({\Bbb A}^2; {\Bbb G}_m^2),
 \ee
where $\PR$ is the Rees bundle construction, and $\Psi$ is the inverse Radon-Penrose transform. 
This is the equivalence claimed in Theorem \ref{main-complex-theorem}. To finish the proof,
it is enough to compare the behavior of the equivalence with respect to the $\Gamma$-action,
see Subsection \ref{complex-version-subsection}.
The definition \eqref{conjugate-structure-equation} of the conjugate Hodge structure implies that
 the coordinates $v_0, v_1, v_2$ on $\check {\Bbb P}^2$ 
associated to the three filtrations in the Rees construction, are transformed under
$\sigma$ 
as follows:
$$v_0\mapsto \overline{v}_0, \quad v_1\mapsto \overline{v}_2,\quad  v_2\mapsto\overline{v}_1.$$
This translates into the action \eqref{action-A2-equation} on the plane ${\Bbb A}^2$ formed by lines
in $\check{\Bbb P}^2$ not meeting $\check{\Bbb P}^0_{F' F''}$. This finishes
the proof of Theorem \ref{main-complex-theorem} and thus of Theorem 
\ref{main-theorem}(a).

\vskip .2cm

To see part (b) of Theorem \ref{main-theorem},
 note that by general properties of the inverse Radon-Penrose
transform, see \eqref{E-lambda-equation},
 and by \eqref{E-gr-equation} , we have
\be\label{E0-fiber-equation} 
E_0 \quad =\quad H^0(\check{\Bbb P}^1_W, \cE) \quad = \quad \gr_\bullet^W(V).
 \ee
Next, to see part (c), note that by the definitions of the Rees module and bundle,
$V$ is recovered as the fiber 
\be 
V\quad  =\quad  \Rs(V, W_\bullet, F'{}^\bullet, F''{}^\bullet)\bigl/ (v_1-1, v_2-1, v_3-1) 
\quad = \quad \cE_{[1:1:1]}.
\ee
The point $[1:1:1]\in\check{\Bbb P}^2$ corresponds to the line 
$$\bigl\{ 1+t_1+t_2=0\bigr\} \quad\subset\quad {\Bbb A}^2,$$
which, in the presence of the real structure $t_1=t, t_2=\overline{t}$,
can be described as $\Ree(t) = -1/2$. Our statement then follows
from the definition of the Radon-Penrose transform,
see Proposition \ref{algebraic-solutions-proposition}.

\vskip .2cm 

To see part (d), notice that flatness of $\nabla$ is equivalent to the property
that $\cE^0$ (and thus $\cE$) is trivial as a vector bundle.
The fact that triviality of the Rees bundle is equivalent to splitting
of the Hodge structure, was pointed out in \cite{P2}, (2.10), Th.2.

\subsection {Noncommutative differential operators and absolute Hodge cohomology.}
Here we prove part (e) of Theorem \ref{main-theorem}. We will prove the following
complex version. Taking into account the real structures is straightforward.

\begin{thm}\label{absolute-Hodge-complex-theorem}
If $(E,\nabla)$ is the ${\Bbb G}_m^2$-equivariant
bundle with connection on ${\Bbb A}^2$ corresponding to
a complex Hodge structure $(V, W_\bullet, F'{}^\bullet, F''{}^\bullet)$, then
we have a natural quasiisomorphism
$$\Gamma\bigl( {\Bbb A}^2, \, \bigl\{ E\buildrel\nabla\over\longrightarrow
 E\otimes\Omega^1_{{\Bbb A}^2}\bigr\}\bigr)
^{{\Bbb G}_m^2} \quad \sim \quad R\Hom_{\MHS_{\Bbb C}} ( {\Bbb C}(0), V).$$
\end{thm}

For the proof we embed $\Bun_\nabla({\Bbb A}^2; {\Bbb G}_m^2)$ into a larger
abelian category in which the trivial bundle $\cO_{{\Bbb A}^2}$ (which corresponds to
the Hodge structure ${\Bbb C}(0)$) has a projective resolution. 

Recall \cite{K} that for any smooth algebraic variety $X/{\Bbb C}$ there is a sheaf
${\Bbb D}_X$ of noncommutative rings on $X$ called the sheaf of
{\it noncommutative differential operators}. We will need the following
properties of ${\Bbb D}_X$. First, ${\Bbb D}_X$ has a 
multiplicative filtration $\{ {\Bbb D}_X^{\leq d}\}$ (by ``order'') with quotients
identified with
\be 
{\Bbb D}_X^{\leq d}/{\Bbb D}_X^{\leq d-1} \quad \simeq\quad T_X^{\otimes \, d}, 
 \ee
 so that the associated graded algebra of $\mathbb D_X$ is the tensor algebra
 of $T_X$.
 
 Second, there is a natural embedding
\be 
\epsilon: T_X\to {\Bbb D}_X^{\leq 1},
 \ee
splitting the first level of the filtration. 

Third, ${\Bbb D}_X$ plays the same role
for  nonflat connections as the ordinary sheaf of differential
operators does for flat ones. To be precise, we have the following.

\begin{prop}\label{nonflat-D-modules-proposition}
 Let $E$ be any quasi-coherent sheaf of $\Oc_X$-modules. Then structures
of a $\DD_X$-module on $E$ extending the $\Oc_X$-module structure, are in
bijection with connections (flat or not) on $E$. \qed
\end{prop}

\begin{exa} Let $X= {\Bbb A}^n$ with coordinates $t_1, ..., t_n$.
Then the ring ${\Bbb D}({\Bbb A}^n)$ is generated by the polynomial ring
${\Bbb C}[t_1, ..., t_n]$ and by the symbols $\nabla_1, ..., \nabla_n$, which are
required to satisfy
$$[\nabla_i, t_j] = \delta_{ij},\quad 1\leq i,j\leq n,$$
and no other relations. In particular, $\nabla_1, ..., \nabla_n$ generate
a free associative algebra. Given a bundle with connection $(E,\nabla)$ on ${\Bbb A}^n$,
the generator $\nabla_i$ acts in $E$ by the covariant derivative $\nabla_{\partial/\partial t_i}$.
\end{exa}

\begin{prop}\label{spencer-prop}
 Let $E'$ be any sheaf of $\DD_X$-modules quasi-coherent over $\Oc_X$
 (i.e., a quasi-coherent $\Oc_X$-module with a connection). Then 
the following 2-term version of
the Spencer sequence is a locally free left ${\Bbb D}_X$-resolution of $E'$:
$$ \begin{aligned}
\Sc^\bullet(E')\,\,=\,\,\bigl\{ {\Bbb D}_X\otimes_{\cO_X} T_X\otimes_{\Oc_X} E' \buildrel d\over\longrightarrow {\Bbb D}_X\bigr\}\otimes_{\Oc_X} E'\}, \\
\quad 
d(P\otimes v\otimes e') = P\cdot \epsilon(v)\otimes e'-P\otimes \epsilon(v)(e').
\end{aligned}
$$
\end{prop}

\noindent {\sl Proof:} Filtering $\Sc^\bullet(E')$  by the subcomplexes
$$\bigl\{ {\Bbb D}_X^{\leq (d-1)}\otimes_{\cO_X} T_X\otimes_{\Oc_X}E' 
\buildrel d\over\longrightarrow 
{\Bbb D}_X^{\leq d}\otimes_{\Oc_X}E'\bigr\},
\quad d\geq 0,$$
we find the associated graded complex to be
\be\label{tensor-algebra-resolution-equation} 
 \biggl( \bigoplus_{d\geq 0} T_X^{\otimes(d-1)}\biggr)\otimes T_X\otimes E' \longrightarrow \bigoplus_{d\geq 0}
T_X^{\otimes d}\otimes E'.
\ee
As the associated graded algebra of $\DD_X$ is the tensor algebra of $T_X$,
we see that the differential in \eqref{tensor-algebra-resolution-equation}
is the tensor multiplication, and so it is
  a resolution of $T_X^{\otimes 0}\otimes E' = E'$ as a left module over the tensor
  algebra.
 Using the exactness of the associated graded complex 
 \eqref{tensor-algebra-resolution-equation}, we deduce the exactness of 
  $\Sc^\bullet(E')$ by a spectral sequence
 argument. \qed

\vskip .3cm

Let now $\Mc$ be the category of $\GG_m^2$-equivariant quasi-coherent sheaves
of $\Oc_{\mathbb A^2}$-modules, equipped with an equivariant connection
(not necessarily flat). In other words, $\Mc$ is the category of $\GG_m^2$-equivariant
$\DD_{\mathbb A^2}$-modules, quasi-coherent over $\Oc_{\mathbb A^2}$. 
Clearly, $\Mc$ is an abelian category. Theorem 
\ref{main-complex-theorem} realizes $\MHS_\CC$ 
as a full subcategory in $\Mc$ (formed by quasi-coherent sheaves which
are vector bundles). 
Note that this subcategory is also abelian, and is closed under extensions in $\Mc$. 
This implies that for any two objects $V, V'\in\MHS_\CC$
with the corresponding objects $E,E'\in\Mc$, the natural morphism
$$\on{Ext}^i_{\MHS_\CC}(V', V) \lra \on{Ext}^i_{\Mc}(E', E)$$
is an isomorphism for $i=0,1$. 

\begin{prop}\label{MHS-M-prop}
 (a) The category $\MHS_\CC$ has cohomological dimension 1, i.e.,
for any two objects $V, V'$ we have $\on{Ext}^i_{\MHS_\CC}(V',V)=0$ for
$i\geq 2$. 

(b) The category $\Mc$ also has cohomological dimension 1.  Therefore for
any two objects $V', V\in\MHS_\CC$ with corresponding
equivariant bundles with connections $E', E$, the natural morphism
$$R\Hom_{\MHS_\CC}(V', V)\lra R\Hom_{\Mc}(E', E)$$
is a quasi-isomorphism.
\end{prop}

\noindent {\sl Proof:} Part (a) 
  was proved by Carlson \cite{Carlson}. Let us prove part (b).
   Given two objects $E, E'$
  of $\Mc$, the complex of vector spaces $R\Hom_{\Mc}(E', E)$
  is obtained from the complex  
  \be\label{RHOM-equation} 
  \ul{R\Hom}_{\DD_{\mathbb A^2}}(E', E)\,\,=\,\,\ul\Hom_{\DD_{\mathbb A^2}}(\Sc^\bullet(E'), E)
  \ee
  of sheaves on $\mathbb A^2$ by taking 
  (the derived functor of) global sections and then
  taking the (derived functor of the) subspace of invariants with respect to $\GG_m^2$.
  Now, the complex \eqref{RHOM-equation} is a 2-term complex of quasi-coherent sheaves
  on $\mathbb A^2$, equivariant with respect to $\GG_m^2$.
  Since $\mathbb A^2$ is affine and $\GG_m^2$ is reductive, we do not
  need to derive the functors of global sections and invariants,
  so $R\Hom(E',E)$ is still a 2-term complex.\qed
  
  \vskip .2cm
  
  This implies Theorem \ref{absolute-Hodge-complex-theorem}.
  Indeed, the statement and the proof of Proposition  \ref{MHS-M-prop} imply that
   $$R\Hom_{\MHS_\CC}(\CC(0), V) = R\Hom_\Mc(\Oc_{\mathbb A^2}, E)
   =\Gamma(\mathbb A^2, \ul\Hom_{\DD_{\mathbb A^2}}(\Sc^\bullet(\Oc_{\mathbb A^2}),
   E))^{\GG_m^2},$$
  which is nothing but the global equivariant de Rham sequence in the formulation of the theorem.
   \qed

 \vskip .2cm

\vfill\eject

\section{The Deligne operator as the holonomy of an
 equivariant connection.}\label{Deligne-operator-section}

\subsection {Reminder on the Deligne operator.}
We start by recalling the description of $\MHS_{\Bbb C}$ by means of the
``Hodge monodromy'' operator $\delta$ given by Deligne \cite{D2}.
For any object $(V, W_\bullet, F'{}^\bullet, F''{}^\bullet)$, there are two splittings
\be\label{two-splittings-equation} 
V\quad =\quad \bigoplus_{p,q\in {\Bbb Z}} V_{F'}^{p,q} \quad = \quad 
\bigoplus_{p,q\in {\Bbb Z}}V_{F''}^{p,q},
\ee
(see \cite{D2}, p. 510) with the following properties. First, both these
splittings induce the same filtration $W_\bullet$ by
\be\label{weight-filtration-equation}
\bigoplus_{p+q\leq n} V_{F'}^{p,q} \quad =\quad \bigoplus_{p+q\leq n} V_{F''}^{p,q} 
\quad = \quad W_n(V).
\ee
Second, the first splitting induces $F'{}^\bullet$ and the second $F''{}^\bullet$
by
\be\label{splittings-filtrations-equation}
F'{}^p(V) = \bigoplus_{p'\geq p, \,\, q\in {\Bbb Z}} V_{F'}^{p,q}, \quad \quad 
F''{}^q (V) = \bigoplus_{p\in {\Bbb Z},\,\, q'\geq q} V_{F''}^{p, q'}.
\ee
Third,
\be 
V_{F'}^{q,p} \equiv \sigma^{-1}(\sigma(V)_{F''}^{p,q})
 \quad {\operatorname{ mod}} \quad W_{p+q-1}(V).
\ee
Here $\sigma$ is the antilinear isomorphism 
\eqref{sigma-map-equation}
 and $\sigma(V)^{p,q}_{F''}$
is the second splitting but for the conjugate Hodge structure,
as defined by \eqref{conjugate-structure-equation}. 
 Let us also denote for simplicity
\be 
\gr^{p,q}(V) = \gr^p_{F'} \gr^q_{F''}\gr_{p+q}^W(V).
 \ee
Then the projections
\be
a_{F'}^{p,q}: V^{p,q}_{F'}\to \gr^{p,q}(V), \quad a_{F''}^{p,q}:  V_{F''}^{p,q}\to
\gr^{p,q}(V)
 \ee
are isomorphisms. So their direct sums, denoted $a_{F'}, a_{F''}$ are isomorphisms 
\be\label{a-F-equation}
 a_{F'}, a_{F''}: V\to \gr^W_\bullet(V) = \bigoplus_{p,q} \gr^{p,q}(W),
 \ee
and so we have the automorphism
\be \delta = a_{F''}a_{F'}^{-1}: \gr_\bullet^W(V)\to \gr^W_\bullet(V),
 \ee
which is known to satisfy
\be (\delta-1)(\gr^{p,q}(V)) \quad \subset \quad \bigoplus_{p'<p, \, q'<q}
\gr^{p', q'}(V).
 \ee
Deligne's characterization is then as follows.

 \begin{thm}\label{deligne-theorem}
  Let $\Delta$ be the category of
pairs $(V^{\bullet \bullet}, \delta)$, where $V^{\bullet\bullet}  = \bigoplus V^{p,q}$
is a finite-dimensional bigraded $\Bbb C$-vector space, and $\delta$ is an
automorphism of $V^{\bullet\bullet}$ satisfying
$$(\delta-1) (V^{p,q}) \quad\subset\quad \bigoplus_{p'<p,\,\, q'<q} V^{p', q'}.$$
Then the functor
$$(V, W_\bullet, F'{}^\bullet, F''{}^\bullet)\,\,\mapsto\,\,
(V^{\bullet\bullet} = \gr^{\bullet\bullet}(V), \,\,\delta =  a_{F''}a_{F'}^{-1})$$
is an equivalence of tensor categories $\MHS_\CC\to\Delta$. 
\end{thm}

\subsection{The Hodge group.} Now, still following \cite{D2}, we reformulate
Theorem \ref{deligne-theorem}
 in terms of representations of appropriate Lie algebras and groups.
Indeed, $\delta$ being unipotent, specifying $\delta$ is equivalent to specfying
its logarithm
\be D= \log (\delta) = \sum_{p,q\geq 1} D_{p,q}, 
 \ee
where $D_{p,q}$ is the bihomogeneous part of  degree $(-p, -q)$, i.e.,
$D_{p,q}(V^{p', q'})\subset V^{p'-p,\,\, q'-q}$ for any $p', q'$. Now, the
operators $D_{p,q}$ can be given arbitrarily subject only to the homogeneity conditions. 
So we get:

 \begin{refo}\label{deligne-group-refo}
 $\MHS_{\Bbb C}$ is equivalent to the category of
finite-dimensional bigraded representations of the free Lie algebra
$${\frak L} = \FLie \bigl\{ z_{p,q}\bigr| \,\, p,q\geq 1 \bigr\},$$
where the generator $z_{p,q}$ has bidegree $(-p, -q)$. 
\end{refo}

The bigrading in $\frak L$ induces the weight filtration $W$ on it, as in 
\eqref{weight-filtration-equation},
and this filtration is compatible with the Lie algebra structure. Each
quotient
 ${\frak L}/W_{-d}{\frak L}$ is finite-dimensional and nilpotent. Denote by
$\exp({\frak L}/W_{-d}{\frak L})$ the corresponding unipotent algebraic group.
We have then the pro-algebraic group
\be {\frak U} \quad = \quad \lim_{\longleftarrow} {} _d \,\,\,\exp({\frak L}/W_{-d}{\frak L}).
 \ee
The bigrading on ${\frak L}/W_{-d}{\frak L}$ can be interpreted as an action of
${\Bbb G}_m^2$, and these actions induce an action on $\frak U$, so we have the semidirect
product
\be{\frak G} = {\Bbb G}_m^2 \ltimes {\frak U}.
\ee
Then,  $\MHS_{\Bbb C}$ is equivalent to the category of finite-dimensional
representations of $\frak G$. The group $\frak G$, sometimes referred to as
the {\it  Hodge group}, is the  group of automorphisms of the fiber functor 
\be\omega_W: \MHS_{\Bbb C}\to \Vect_{\Bbb C}, \quad (V, W_\bullet, F'{}^\bullet, F''{}^\bullet)
\mapsto  \gr^{\bullet\bullet}(V).
\ee

 \subsection {Comparison of two descriptions.}
  We now compare the description
of $\MHS_{\Bbb C}$ given by Theorem \ref{deligne-theorem} (or, equivalently, by 
Reformulation \ref{deligne-group-refo})
with the description of Theorem \ref{main-complex-theorem} via equivariant connections.

 \begin{thm}\label{delta-holonomy-theorem}
  Let $(V, W_\bullet, F'{}^\bullet, F''{}^\bullet)\in\MHS_{\Bbb C}$,
and $(E,\nabla)$ be the corresponding ${\Bbb G}_m^2$-equivariant bundle with connection
in ${\Bbb A}^2$. Then, under the identification $E_0=\gr^{\bullet\bullet}(V)$ given by
\eqref{E0-fiber-equation}, the Deligne operator
 $$\delta: E_{(0,0)} = V^{\bullet\bullet}\to V^{\bullet\bullet} = E_{(0,0)}$$
 is recovered
as the holonomy of $\nabla$ along the boundary of the triangle with vertices 
$(0,0)$, $(-1,0)$, $(0, -1)$.
\end{thm}

\noindent {\sl Proof:} 
We realize the bundle $\cE$ explicitly by a patching function. 
We represent the punctured projective plane
 as the union of two affine
planes:
\be 
\check{\Bbb P}^2_0 \quad  =\quad  \check{\Bbb P}^2 - \check{\Bbb P}^0_{F'F''} \quad = \quad
\check{\Bbb A}^2_{F'} \cup \check {\Bbb A}^2_{F''},
\ee
where 
\be
\begin{aligned}
\check{\Bbb A}^2_{F'}  = \Spec \, {\Bbb C}[\xi_0, \xi_1], \quad \xi_0 = v_0/v_2, \quad \xi_1 = v_1/v_2, \\
\check{\Bbb A}^2_{F''} = \Spec \, {\Bbb C}[ \eta_0, \eta_2], \quad \eta_0 = v_0/v_1,
\quad  \eta_2= v_2/v_1,
\end{aligned}
\ee
are the affine charts centered at the points $\check{\Bbb P}^0_{WF'}$ and
$\check{\Bbb P}^0_{WF''}$. We are interested only in the restriction $\cE^0$ of $\cE$
to $\check{\Bbb P}^2_0$, as points of ${\Bbb A}^2$ correspond to lines contained
in $\check{\Bbb P}^2_0$. 

The restriction of the Rees bundle $\cE$ to $\check{\Bbb A}^2_{F'}$ is just the
Rees bundle on ${\Bbb A}^2$ corresponding to the pair of filtrations $W,F'$,
see \cite{P2} (2.6.1). In other words, this is the bundle correponding to the
 ${\Bbb C}[\xi_0, \xi_1]$-module
\be 
M_{F'} \quad = \quad \bigoplus_{i,j} (W_i\cap F'{}^{-j})(V)\,\, \xi_0^i \xi_1^j 
\quad\subset\quad V[\xi_0^{\pm 1}, \xi_1^{\pm 1}]. 
\ee
Note that $W_\bullet, F^\bullet$ are simultaneously split by the first  bigrading
$(V^{p,q}_{F'})$ from 
\eqref{two-splittings-equation}, so by 
\eqref{splittings-filtrations-equation} we can identify $M_{F'}$ with a
free module as follows:
\be\label{M-F'-equation}
M_{F'}\quad  =\quad  \bigoplus_{i,j} \bigoplus_{p+q\leq i \atop p\geq -j} V_{F'}^{p,q}\,\, \xi_0^i \xi_1^j 
\quad  =\quad  
\bigoplus_{p,q} V_{F'}^{p,q}\otimes\,\,  \xi_0^{p+q} \xi_1^{-p}\, {\Bbb C}[\xi_0, \xi_1].
 \ee
In other words, we have the  trivialization
\be\label{tau-F'-equation}
 \begin{split}
\tau_{F'}: \gr^{\bullet\bullet}(V) \otimes {\Bbb C}[\xi_0, \xi_1] \lra M_{F'}, \hskip 3cm\\
\tau_{F'}(v_{pq}\otimes f(\xi_0, \xi_1)) = a_{F'}^{-1}(v_{pq})\otimes \xi_0^{p+q}\xi_1^{-p} f(\xi_0, \xi_1),
\quad v_{p,q}\in\gr^{p,q}(V).
\end{split}
\ee
Here $a_{F'}^{-1}$ is the inverse of the isomorphism induced by the splitting $(V^{p,q}_{F'})$, see 
\eqref{a-F-equation}.

Similarly, the restriction of $\cE$ to $\check{\Bbb A}^2_{F''}$ is the Rees
bundle on ${\Bbb A}^2$ corresponding to $W, F''$, so it corresponds to the
${\Bbb C}[\eta_0, \eta_2]$-module
\be\label{M-F''-equation} 
M_{F''}\quad  =\quad  \bigoplus_{k,l} (W_k\cap F''{}^{-l})(V)\,\, \eta_0^k\eta_2^l\quad  =\quad
\bigoplus_{k,l} \bigoplus_{p+q\leq k,\atop q\geq -l} V^{p,q}_{\barF} \eta_0^k \eta_2^l,
 \ee
where we now used the second splitting in 
\eqref{two-splittings-equation}.  So we get a trivialization
\be\label{tau-F''-equation}
 \begin{split}
\tau_{F''}: \gr^{\bullet\bullet}(V)\otimes {\Bbb C}[\eta_0, \eta_2] \to M_{F''},
\hskip 2cm
\\
\tau_{F''}(v_{pq}\otimes g(\eta_0, \eta_2)) = a_{F''}^{-1}(v_{p,q})\otimes \eta_0^{p+q}\eta_2^{-q} g(\eta_0,
\eta_2).
\end{split}
\ee
This means that we have glued $\cE^0$ out of two trivial bundles, each  with fiber $\gr^{\bullet\bullet}(V)$. 
Let us view a point $\xi = (\xi_0, \xi_1)\in \check{\Bbb A}^2_{F'}$ with both coordinates
nonzero as an element of the torus ${\Bbb G}_m^2$. Let $\xi\mapsto \phi(\xi)$ be the
torus action on $\gr^{\bullet\bullet}(V)$ corresponding to the bigrading.
In coordinates $\xi_0, \xi_1$ it has the form:
\be 
\phi(\xi) (v_{p,q}) = \xi_0^{p+q} \xi_1^{-p} v_{p,q}, \quad v_{p,q}\in \gr^{p,q}(V).
 \ee
Writing now the identification of these two bundles on $\check{\Bbb A}^2_{F'}\cap\check{\Bbb A}^2_{F''}$,
i.e.,  expressing the $\eta$'s through the $\xi$'s and accounting for the monomial
factors in \eqref{tau-F'-equation} and \eqref{tau-F''-equation},
we find straightforwardly:

 \begin{prop}\label{patching-function-proposition}
 The patching function  $\Phi = \tau_{\barF}^{-1} \tau_F$
of $\cE^0$ with respect to the above
trivializations has the form
$$\Phi(\xi) = \phi(\xi)^{-1} \, \delta \, \phi(\xi).\qed$$
\end{prop}

We now prove Theorem \ref{delta-holonomy-theorem}.
 Recall that  each point $T=(t_1, t_2)\in {\Bbb A}^2$
corresponds to a line
\be 
\lambda_T = \bigl\{ v_0+t_1 v_1 + t_2 v_1 =0\bigr\} \quad\subset\quad
\check{\Bbb P}^2_0.
\ee
Given two distinct points $T, T'$, we denote by $[T,T']$ the
straight line segment joining them. Since $\cE$ is trivial on each
$\lambda_T$, the restriction map
\be 
R_{TT'} : E_T = H^0(\lambda_T, \cE) \longrightarrow \cE_{\lambda_T \cap\lambda_{T'}} =
H^0_\nabla([T, T'], E),
\ee
is an isomorphism. By general properties of the Radon-Penrose transform
(cf, e.g., \cite{WW}, p. 377), the holonomy along $[T, T']$ is found as
\be H_{TT'} = R^{-1}_{T'T}\circ R_{TT'} : E_T \longrightarrow E_{T'}.
\ee
Now notice that $\lambda_{(0,0)} = \check{\Bbb P}^1_W$, while
$\lambda_{(-1,0)}$ is the line $\{ v_0=v_1\}$ joining the points $\check{\Bbb P}^0_{WF}
=[0:0:1]$ and $[1:1:1]$. Similarly, $\lambda_{(0, -1)}$ is the line 
$\{v_0=v_2\}$ joining $\check{\Bbb P}^0_{W\barF} = [0:0:1]$ and $[1:1:1]$. 
By Proposition \ref{patching-function-proposition},
 the value of the patching function $\Phi(\xi)$ at the point
$\xi_0=\xi_1=1$, which is the same as $[1:1:1]$, is equal to $\delta$. 
This implies:

 \begin{lem} If we identify $H^0(\lambda_{(-1,0)}, \cE)$ and
$H^0(\cE, \lambda_{(0, -1)}, \cE)$ 
with $\gr^{\bullet\bullet}(V)$, using the trivializations
$\tau_{F'}$ and $\tau_{F''}$, then the composite isomorphism
$$H^0(\lambda_{(-1,0)}, \cE)\buildrel R_{(-1,0), (0, -1)} \over\longrightarrow
\cE_{[1:1:1]} \buildrel R^{-1}_{(0, -1), (-1, 0)} \over\longrightarrow
H^0(\lambda_{(0, -1)}, \cE),$$
is identified with $\delta$.
\end{lem}

The lemma implies Theorem \ref{delta-holonomy-theorem},
 as with respect to our identifications the
holonomies along the two other sides of the triangle are equal to 1.

\subsection { Equivariant connections in coordinates.} Let
$V^{\bullet\bullet}$ be a finite-dimensional
bigraded $\Bbb C$-vector space, and $\widetilde{V}^{\bullet\bullet}$
be the corresponding ${\Bbb G}_m^2$- equivariant vector bundle on ${\Bbb A}^2$. 
As a vector bundle, it is trivial: $\widetilde{V}^{\bullet\bullet}= V^{\bullet\bullet}
\otimes\cO_{{\Bbb A}^2}$, while the ${\Bbb G}_m^2$-action is given by
\be (\lambda_1, \lambda_2) \bigl( v_{pq}\otimes f(t_1, t_2)\bigr) \quad = \quad
t_1^p t_2^q v_{pq}\otimes f(\lambda_1 t_1, \lambda_2 t_2), \quad v_{pq} \in V^{p,q}.
 \ee

\begin{lem} Every ${\Bbb G}_m^2$-equivariant vector bundle on ${\Bbb A}^2$
is equivariantly isomorphic to a bundle of the form $\widetilde{V}^{\bullet\bullet}$. 
\end{lem}

\noindent {\sl Proof:} It is well known that every 
${\Bbb G}_m^2$-equivariant vector bundle on ${\Bbb A}^2$ is the Rees bundle
corresponding to two filtrations, see, e.g., \cite{P2} (2.5.3) Prop. 14. 
Our statement follows from the
fact that any two filtrations can be simultaneously split by a bigrading, see, e.g.,
\cite{P2} (2.1.3), Lemme 1. \qed

 \begin{prop}\label{invariant-connections-coordinates-prop}
Any  $\GG_m^2$-invariant connection in  $\widetilde{V}^{\bullet\bullet}$ has the form
$$\nabla = d+\Omega, \quad \Omega = \sum_{p,q\geq 1} A_{p,q} t_1^{p-1} t_2^q dt_1 + B_{p,q}
t_1^p t_2^{q-1} dt_2,$$
where $d$ is the standard flat connection of the trivial bundle $\widetilde{V}^{\bullet\bullet}$,
and $A_{p,q}, B_{p,q}$ are endomorphisms of ${V}^{\bullet\bullet}$
of bidegree $(-p, -q)$. \qed
\end{prop}

\noindent {\sl Proof:} The space of
$d$-covariantly constant (i.e., constant) sections of 
$\widetilde{V}^{\bullet\bullet}$ is preserved by the action $\GG_m^2$
(although individual elements of this space may not be). Next, the datum of any connection
 at a given point $x\in\mathbb A^2$ is given by the subspace $C_x$ in the space
of germs of sections near $x$ which are covariantly constant up to the first order
of tangency. Looking at the trivial connection $d$, we see,
just as for global sections, that the action of $\lambda\in\GG_m^2$ 
takes the subspace $C_x$ into the subspace $C_{\lambda(x)}$. 
This means that $d$ is in fact a $\GG_m^2$-invariant
connection in $\widetilde{V}^{\bullet\bullet}$. Therefore any other invariant
connection has the form $d+\Omega$, where $\Omega$ is a global 
1-form on $\mathbb A^2$ with values in $\End(\widetilde{V}^{\bullet\bullet})$, which
is $\GG_m^2$-invariant, i.e., has total degree 0. The sum in the proposition
is nothing but the general shape of such a 1-form. \qed

\vskip .3cm

We denote 
\be W_n(V^{\bullet\bullet}) = \bigoplus_{p+q\leq n} V^{p,q},\quad
W_n(\widetilde{V}^{\bullet\bullet}) = \widetilde{W_n(V^{\bullet\bullet})}
\ee
the weight filtration of the bigraded space $V^{\bullet\bullet}$ and of the associated
bundle. Then each connection in Proposition
\ref{invariant-connections-coordinates-prop}
 preserves $W_n(\widetilde{V}^{\bullet\bullet})$
and induces the trivial connection on the quotients. 

Isomorphisms among  connections in Proposition
\ref{invariant-connections-coordinates-prop} which induce the identity
on the quotients of the weight filtration  correspond to gauge transformations
\be\label{gauge-transform-coordinates-equation} 
\Omega\mapsto g^{-1} dg + g^{-1} \Omega g, \quad g = g(t_1, t_2) = \sum_{p,q \geq 0} 
C_{p,q} t_1^p t_2^q,\quad C_{0,0}=1,  
\ee
where $C_{p,q}$ is an endomorphism of  ${V}^{\bullet\bullet}$
of bidegree $(-p, -q)$. 

\begin{prop}\label{Apq+Bpq=0-prop}
 Each equivalence class of connections as in
Proposition \ref{invariant-connections-coordinates-prop}
with respect to transformations \eqref{gauge-transform-coordinates-equation} 
 contains a unique connection
satisfying
$$A_{p,q} + B_{p,q} = 0, \quad \forall p,q.$$
\end{prop}

The proof is easy, by induction on the length of the weight filtration.\qed

\begin{rem} The condition in Proposition 
\ref{Apq+Bpq=0-prop} is
a particular case of the so-called Fock-Schwinger gauge condition in physics
 which for a connection
$$\nabla = d+\sum \cA_\nu(t_1, ..., t_n) dt_\nu$$
 in a trivial bundle over ${\Bbb R}^n$ or ${\Bbb C}^n$, reads:
$$\sum _\nu t_\nu \cA_\nu(t_1, ..., t_n) \quad = \quad 0. $$
\end{rem}

 \subsection {A different set of generators of the Lie algebra $\frak L$. }
We now have two ways of describing a complex  mixed Hodge structure
$(V, W_\bullet, F'{}^\bullet, F''{}^\bullet)$ with fixed
weight quotients $\gr^W_\bullet(V) = V^{\bullet\bullet}$.
The first one is by means of the Deligne operator $\delta$ or,
equivalently, of the components $D_{p,q}, \, p,q\geq 1$, of
$D=\log (\delta)$. The second one is by means of an equivariant connection, which,
by Proposition \ref{Apq+Bpq=0-prop}
 we can uniquely represent by the gauge potential
\be 
\Omega =  \sum_{p,q\geq 1} A_{p,q} (t_1^{p-1}t_2^q dt_1 - t_1^p t_2^{q-1} dt_2),
 \ee
where the $A_{p,q}$ satisfy the same homogeneity conditions as the $D_{p,q}$
and, apart from these conditons, can be taken arbitrarily. 
The connection $\nabla = d+\Omega$ is trivial along the coordinate axes $t_i=0$, so
Theorem \ref{delta-holonomy-theorem}
 gives the following relation between the $D_{p,q}$ and the $A_{p,q}$:
\be 
D = \sum_{p,q\geq 1}  D_{p,q} = \log P\exp \int_{(-1,0)}^{(0, -1)} \Omega.
\ee
Here we used the notation $P\exp\int$ (path-ordered exponential)
 for the holonomy of a connection. 
This implies the following.

 \begin{prop} There exist universal relations
$$D_{p,q}\quad  =\quad  (-1)^{p+q} 
{p+q\choose p} A_{p,q} + S_{p,q} \bigl(\bigl\{ A_{p', q'}\bigl| \, p'<p,\, q'<q
\bigr\}\bigr),$$
where $S_{p,q}$ are Lie polynomials with rational coefficients
in the lower $A_{p',q'}$,  bihomogeneous
of bidegree $(-p, -q)$. \end{prop}

To prove the proposition, we repeat the above reasoning in the
``universal''  situation.
 Consider the free bigraded Lie algebra
 \be {\frak L}' = \FLie \bigl\{ \alpha_{p,q}\bigl| \,\, p,q\geq 1, \,\,
\deg (\alpha_{p,q}) = (-p, -q)\bigr\},
\ee
similarly to the one in Reformulation \ref{deligne-group-refo},
and let ${\frak L}'_{\leq n}$ be the subalgebra generated by the $\alpha_{p,q}$ with
$p+q\leq n$. Let $R_n$ be the completed universal enveloping algebra of ${\frak L}'_{\leq n}$,
i.e., the algebra of noncommutative formal power series in the $\alpha_{p,q}$,
$p+q\leq n$. We have then a connection on ${\Bbb A}^2$ with values in $R_n$:
\be 
\nabla = d+\omega, \quad \omega = \sum_{p+q\leq n} \alpha_{p,q} 
(t_1^{p-1}t_2^q dt_1 - t_1^p t_2^{q-1} dt_2).
 \ee
The holonomy of this connection along the segment $[(-1,0), (0, -1)]$ is
then a well defined element of $R_n$, and  we consider its logarithm $z$
and the bihomogenous components $z_{p,q}$ of $z$:
$$z = \sum_{p,q\leq n} z_{p,q} = \log P\exp\int_{(-1,0)} ^{(0, -1)} \omega \leqno (3.5.6)$$
Since $z$ is a primitive element, it is a Lie series in the 
 $\alpha_{p,q}, p+q\leq n$. By  degree considerations, each $z_{p,q}$
is a Lie polynomial.  
Now, modulo the commutators, there is no difference between the path ordered
exponential and the usual exponential of the integral, so we have
$$ z 
\quad \equiv \quad \int_{(-1,0)} ^{(0, -1)} \omega \quad = \quad 
\sum_{p+q\leq n} \alpha_{p,q} \int_{-1}^0 -t^{p-1} (-1-t)^{q-1} dt,$$
whence the claim. \qed

\vskip .2cm

We can view the variables $z_{p,q}$ as the generators of the free Lie algebra
$\frak L$ from Reformulation \ref{deligne-group-refo}.
 So taking the logarithm of the holonomy defines an
homomorphism
\be 
{\frak L}\to {\frak L'}, \quad z_{p,q} \mapsto (-1)^{p+q} {p+q\choose p}\alpha_{p,q} + S_{p,q}
\bigl(\bigl\{\alpha_{p', q'}\bigl| \, p'<p,\, q'<q
\bigr\}\bigr),
\ee
which, because of its triangular form, is an isomorphism. In particular,
the  $\alpha_{p,q}$ can be expressed back through the $z_{p,q}$ and
provide an alternative system of bihomogeneous generators of $\frak L$.

 \begin{rem}\label{twistor-line-remark}
  The generators $\alpha_{p,q}$
essentially  coincide with
the generators introduced by Goncharov \cite{G}
starting from totally different principles. In fact, the main
feature in Goncharov's  approach to mixed Hodge structures  is a connection
on the  affine line ${\Bbb A}^1$  (called  the ``twistor line'' in \cite{G}) 
whose holonomy along $[0,1]$ is  equal to the Deligne operator $\delta$.
In our approach, this connection appears as the restriction of
a ${\Bbb G}_m^2$-equivariant connection from ${\Bbb A}^2$ to the
affine line $t_1+t_2=-1$. 
 \end{rem}

\subsection {Geometric interpretation of $\frak L$.}
 We would now like to give a
different interpretation of the bigraded free Lie algebra $\frak L$.
 Namely, we observe that
$\frak L$ {\it can be  identified with the commutant of the free Lie algebra 
on two generators}.

\vskip .2cm

To be precise, let $A$ be a $\Bbb C$-vector space, and
\be \FLie (A) = \bigoplus_{d\geq 1} \FLie_d(A)
 \ee
be the free Lie algebra on $A$ graded by the degree of commutator monomials,
i.e., by putting $A$ in degree 1 and requiring that the bracket be
compatible with the grading. Assume that $A$ is finite-dimensional.
Then the algebraic group $GL(A)$ acts on
$\FLie(A)$ by Lie algebra automorphisms., and we are interested in the
action on the commutant
\be 
[\FLie(A),\, \FLie(A)] = \FLie_{\geq 2}(A).
\ee
This commutant, considered as an abstract Lie algebra is free.
Indeed, by the Shirshov-Witt theorem (\cite{R}, Ch.2)  any subalgebra of
a free Lie algebra is free. However, there is no canonical choice
of free generators for the commutant. For different choices the
``spaces of generators'' are identified with each other via the
identification with the first homology space (``indecomposable elements'')
\be\label{first-homology-equation}
H_1^{\Lie}(\FLie_{\geq 2}(A), {\Bbb C}) = \FLie_{\geq 2}(A)/[\FLie_{\geq 2}(A),
\FLie_{\geq 2}(A)]. 
\ee
In other words, if $B\subset \FLie_{\geq 2}(A)$ is a graded subspace, then
the following are equivalent:

\vskip .2cm

\begin{enumerate}
 \item The natural map $\FLie(B)\to \FLie_{\geq 2}(A)$ is
an isomorphism, so $B$ is a space of free generators;

\vskip .1cm

 \item The projection $B\to H_1^{\Lie}(\FLie_{\geq 2}(A), {\Bbb C})$
is an isomorphism.

\end{enumerate}

\vskip .2cm

The first homology space \eqref{first-homology-equation}
 was identified, as a $GL(A)$-module, by
Reutenauer (\cite{R}, (8.6.12)). His result reads:
\be 
H_1^{\Lie}(\FLie_{\geq 2}(A), {\Bbb C}) = \bigoplus_{d\geq 1} \Sigma^{d,1}(A),
 \ee
where $\Sigma^{d,1}$ is the Schur functor (irreducible representation of
$GL$) corresponding to the Young diagram $(d,1)$. It is also well known
(see, e.g., \cite{GKZ}, Prop. 14.2.2) that as a $GL(V)$-module, 
\be 
\bigoplus_{d\geq 1}\Sigma^{d,1}(A) = \Omega^{2, \cl}_{\pol}(A^*), 
\ee
 is just the space of closed polynomial 2-forms on the affine space $A^*$.

\vskip .2cm

Now, assume that $A = {\Bbb C}^2$ is 2-dimensional, and $t_1, t_2$ form the
standard  basis of $A$. We can consider $t_i$ as a linear function
on $A^*$. 
 Then all 2-forms are closed, so a basis of
$\Omega^{2, \cl}_{\pol}(A^*)$ is formed by the monomials
\be 
w_{p,q} = t_1^{p-1} t_2^{q-1} dt_1 dt_2, \quad p,q\geq 1. 
\ee
If we equip $\FLie({\Bbb C}^2)$ with the bigrading starting with
\be 
{\deg}(t_1) = (-1,0), \quad \deg (t_2) = (0, -1),
\ee
and then restrict this bigrading to $\FLie_{\geq 2}({\Bbb C}^2)$, 
then $w_{p,q}$ has the same bidegree as the generator $z_{p,q}\in{\frak L}$,
namely $(-p, -q)$. So there is  a graded isomorphism  of Lie algebras
\be 
\phi: {\frak L}\to \FLie_{\geq 2}({\Bbb C}^2).
\ee
There is, however, a choice in constructing such a $\phi$, 
which is a choice of lifting of  each $w_{p,q}$
to an element $\phi(z_{p,q})$ of $\FLie_{\geq 2}({\Bbb C}^2)$. One possible
way to fix these  liftings is  by
putting:
\be 
\phi(z_{p,q}) =  [t_2, [t_2, ..., [t_2, [t_1, [t_1, ..., [t_1, t_2]...]=
\ad(t_2)^{q-1}(\ad(t_1)^p(t_2)).
 \ee

\subsection{ $\frak L$ as a  fundamental Lie 
algebra.}\label{fundamental-algebra-subsection} 
It was further shown in (\cite{K},  (4.3)),  that for any finite-dimensional $A$ 
as before,  the Lie algebra
$\FLie_{\geq 2}(A)$ is acted upon not just by $GL(A)$ but by the group
of formal changes of coordinates in $A$, i.e., by
\be 
\Aut \, {\Bbb C}[[s_1, ..., s_n]], \quad {\operatorname{ if}} \quad A^*= \bigoplus_{i=1}^n \, {\Bbb C}\cdot s_i.
\ee
This group is pro-unipotent, while $\FLie_{\geq 2}(A)$  is a
discrete, infinitely generated Lie algebra. The action, when restricted to any
finitely generated subalgebra in $\FLie_{\geq 2}(A)$, factors through
a finite-dimensional unipotent quotient.

Alternatively, this means that to any smooth $n$-dimensional algebraic variety $X/{\Bbb C}$
 and any point $x\in X$
there is a naturally associated {\it fundamental Lie algebra} $\cP(X,x)$
which is isomorphic to $\FLie_{\geq 2}({\Bbb C}^n)$ but not canonically.
Any choice of a formal coordinate system near $x$ gives rise to an isomorphism between
the two. It was also shown that the first cohomology of $\cP(X,x)$ is naturally
identified as follows:
\be\label{fundamental-algebra-equation} 
H^1_{\Lie}(\cP(X,x), {\Bbb C}) = \widehat{\Omega}^{2, \cl}_{X,x},
\ee
where the RHS is the space of formal germs of closed 2-forms on $X$ near $x$. 

Further, there is a  ``nonabelian'' version of 
\eqref{fundamental-algebra-equation}  which  was proved in \cite{K}. 
Let  $\frak g$ be any complex  Lie algebra. Then   $H^1_{\Lie}({\frak g}, {\Bbb C})$
is the same as the group of isomorphism classes of 1-dimensional representations
of $\frak g$, with the operation induced by tensor product. So we denote by
$\Rep({\frak g})$ the tensor category of finite-dimensional representations of $\frak g$.
Then we have an equivalence of tensor categories:
\be 
\Rep(\cP(X,x)) \simeq \widehat{\Bun}_\nabla (X,x), 
\ee
where the RHS is the category of formal germs of vector bundles with
connections on $X$ near $x$, see \cite{K}, Theorem 4.4.3. 
This equivalence is natural with respect to maps of manifolds preserving base
points, in particular, to any algebraic group $G$ acting on $(X,x)$.

\vskip .2cm

We now specialize to $X={\Bbb A}^2$, with $G= {\Bbb G}_m^2$ acting diagonally. We get that
the category of bigraded representations of $\cP({\Bbb A}^2, 0) = {\frak L}$ is
equivalent to the category of ${\Bbb G}_m^2$-equivariant objects in 
$\widehat{\Bun}_\nabla({\Bbb A}^2, 0)$. Now, an equivariant formal germ
of a connection on ${\Bbb A}^2$ near 0 must be  a polynomial one by degree reasons.
So the latter category is identified with $\Bun_\nabla({\Bbb A}^2; {\Bbb G}_m^2)$.
This provides an alternative way to relate our description of
$\MHS_\CC$ with Theorem \ref{deligne-theorem}.

\end{document}